\documentclass[11pt]{article}
\date{}
\newcounter{list}
\newtheorem{Th}{Theorem}[section]
\newtheorem{Lem}[Th]{Lemma}
\newtheorem{Prop}[Th]{Proposition}
\newtheorem{Cor}[Th]{Corollary}
\newtheorem{Ex}[Th]{Example}
\newtheorem{Rem}[Th]{Remark}
\newtheorem{Def}[Th]{Definition}
\newtheorem{Def-Lem}[Th]{Definition and Lemma}
\newtheorem{Def-Th}[Th]{Definition and Theorem}
\newenvironment{listing}[1]{\setcounter{list}{0}
\begin{list}{\rm{#1\arabic{list})}}{\usecounter{list}}}
{\end {list}}
\def\A{\rightarrow}
\def\bar{\underline}
\def\Dot{^{\bullet}}
\def\H#1#2{\mathop{{\rm H}^{#1}({#2}\Dot)}\nolimits}

\def\dif#1#2{\mathop{ { d_{#1}}^{#2}}\nolimits}
\def\ex#1#2#3{0\A{#1}\A{#2}\A{#3}\A 0}
\def\hom#1{\mathop{\rm Hom _{#1}}\nolimits}
\def\hombar#1{\mathop{\bar{\rm Hom} _{#1}}\nolimits}
\def\ext#1#2{\mathop{\rm Ext ^{#1}_{#2}}\nolimits}

\def\mapdown#1{\Big\downarrow \rlap{$\vcenter{\hbox{$\scriptstyle#1$}}$}}

\def\stcong{\stackrel{\scriptstyle st}{\cong}}
\def\syz#1#2#3{\Omega ^{#1} _{#2} (#3)}

\def\tr{\mbox{tr}\,}

\def\plim{\mathop{\raise-.2em\hbox{$\def\arraystretch{0}\begin{array}{c}\lim\\
\leftarrow\end{array}$}}}

\def\cok{\mathop{\rm cok}\nolimits}
\def\proj{\mathop{\rm proj}\nolimits}

\def\ll{\mathop{\cal L}\nolimits}

\def\hpy{\mathop{\sf  K}\nolimits}

\def\mod{\mathop{\rm mod}\nolimits}

\def\cok{\mathop{\rm Cok}\nolimits}

\def\ker{\mathop{\rm Ker}\nolimits}
\def\Im{\mathop{\rm Im}\nolimits}

\def\pd{\mathop{\rm pd}\nolimits}

\def\id{\mathop{\rm id}\nolimits}

\def\ass{\mathop{\rm Ass}\nolimits}

\def\tr{\mathop{\rm Tr}\nolimits}
\def\spep#1{\mathop{{}^{\bullet}\strut\kern-.1em{#1}}\nolimits}

\begin {document}
\begin{center}
{\large\bf Morphisms represented by monomorphisms \\ }
{ Kiriko Kato \\
Department of Applied Mathematics, Osaka Women's University\\
Sakai, Osaka 590-0035, JAPAN   \\   
e-mail: { kiriko@appmath.osaka-wu.ac.jp}}
\end{center}

\begin{center}
{\bf Abstract} \end{center}

{\small Every homomorphism of modules is projective-stably equivalent to 
an epimorphism but is not always to a monomorphism. 
We prove that a map is projective-stably equivalent to a monomorphism 
if and only if its kernel is torsionless, that is, a first syzygy. 
If it occurs although, there can be various monomorphisms that are 
projective-stably equivalent to a given map. But in this case 
there uniquely exists a "perfect" monomorphism 
to which a given map is projective-stably equivalent. }

\vspace{5mm}

2000 Mathematics Subjects Classification: 13D02, 13D25, 16D90

\section{Introduction} 
Let $R$ be a semiperfect ring. 
A morphism $f:A \A B$ and $f': A' \to B'$ in $\mod R$ are said to be 
projective-stably equivalent if they are isomorphic in $\bar{\mod R}$; 
if there exist morphisms $\alpha :A\A A'$  and 
$\beta :B \A B'$ such that $\bar\alpha$ and $\bar\beta$ are isomorphisms 
and $\bar{\beta \circ f} = \bar{ f' \circ \alpha}$ in $\bar{\mod R}$. 
We say a morphism $f$ is represented by monomorphisms ("rbm" for short) 
if there exists a monomorphism in $\mod R$ that is projective-stably equivalent to $f$. 

For any homomorphism $f: A\A B$ of $R$-modules, 
$ \pmatrix{f&\rho _B\cr} : A\oplus P_B \A B$ is surjective 
with a projective cover $ \rho _B : P_B \to B$. 
Thus every morphism is represented by epimorphisms. 
The choice of epimorphism is unique; if an epimorphism $f'$ is 
projective-stably equivalent to $f$, 
then $f'$ is isomorphic to $(f ~\rho _B )$ 
up to direct sum of projective modules. 

On the other hand, every morphism is not always represented by monomorphisms. 
Even if a morphism $f$ is rbm, 
the choice of monomorphism is not unique; 
there would be two monomorphisms that are not isomorphic up to direct sum of 
projective modules and both of which are projective-stably equivalent to $f$. 
  
The purpose of the paper 
is finding a condition of a given map to be rbm. 
The problem was posed by Auslander and Bridger \cite{ABr}. 
They proved that a map is rbm if and only if 
it is projective-stably equivalent to a "perfect" monomorphism. 
 An exact sequence 
of $R$-modules is called perfect if its $R$-dual 
is also exact. A perfect monomorphism refers to 
a monomorphism whose $R$-dual is an epimorphism. 
This is our first focal point. 
We studied the situation where a map is rbm, 
especially the structure of monomorphisms into which 
a given map is modified. 
And we obtained the obstruction for a given map 
to be rbm. 
In the case that a map is rbm, 
the choice of a monomorphism is not unique, but then 
a perfect monomorphism projective-stably equivalent to the given map 
is uniquely determined  
up to direct sum of projective modules. (Theorem~\ref{p357}.) 

Our next focus is an analogy to the homotopy category 
$\hpy ( {\mod R} )$ of $R$-complexes. 
In \cite{Okayama} Theorem 2.6, the author showed 
a category equivalence between $\bar{\mod R}$ 
and a subcategory of $\hpy ( {\mod R} )$. 
Due to this equivalence, we describe the obstruction of 
being rbm 
with a homology of a complex associated to the given map. 

Looking at Theorem ~\ref{p357}, we see that 
when a morphism is rbm, 
its pseudo-kernel is always the first syzygy of 
its pseudo-cokernel. So it is tempting to ask if 
torsionlessness of the kernel is equivalent to 
rbm condition. 
This is our third point. Actually, for this we need Gorensteinness. 
\begin{quote} Theorem ~\ref{p284} : 
Suppose the total ring of fractions $Q(R)$ of a ring $R$ is Gorenstein.  
A morphism $f$ is rbm 
if and only if $\ker f$ is torsionless, equivalently, 
a first syzygy. 
\end{quote}

\section{Preliminaries}

We shall fix the notations and give some review on the correspondence between stable module category and homotopy class category 
of complexes. We omit the proofs for results that are 
in \cite{Okayama}.  

Throughout the paper, $R$ is a commutative semiperfect ring, equivalently a finite direct sum of local rings; 
that is, each finite module has a projective cover ( see \cite{Miyachi} 
for semiperfect rings). 
The category of finitely generated $R$-modules is denoted by $\mod R$, 
and the category of finite projective $R$-modules is denoted by $\proj R$.  
By an $R$-module we mean "a finitely generated $R$-module". 
For an $R$-module $M$, $\rho _M : P_M \to M$ 
denotes a projective cover of $M$. 
For an abelian category ${\cal A}$, $\hpy ( {\cal A}) $ stands for 
the category of the homotopy equivalence class of complexes in ${\cal A}$. 
A complex is denoted as  
\[ F\Dot : \cdots \A F^{n-1}  \stackrel{\dif{F}{n-1}}{\A} F^{n}  \stackrel{\dif{F}{n}}{\A} F^{n+1} 
\A \cdots .  \] 
A morphism in ${\hpy ( {\cal A})}$ is a homotopy equivalence class of 
chain maps.  
A trivial complex is a split exact sequences of 
projective modules. 
Truncations of a complex $F\Dot$ are defined as follows:
\[   \tau _{\leq n} F\Dot : \cdots \A F^{n-2}  \stackrel{\dif{F}{n-2}}{\A} F^{n-1} 
\stackrel{\dif{F}{n-1}}{\A} F^n \A 0\A 0 \A \cdots , \]
\[   \tau _{\geq n} F\Dot : \cdots \A 0\A 0\A F^{n}  \stackrel{\dif{F}{n}}{\A} 
F^{n+1} \stackrel{\dif{F}{n+1}}{\A} F^{n+2}  \A  \cdots  \]
An $R$-dual $F^* _{\bullet}$ of a complex $F\Dot$ is the cocomplex 
such as $F^* _n = {(F_n ) }^*$,  
$ {\mathop{d}}^{F^*}_n = { (\dif{F}{n-1})}^*$ where ${}^*
$ means $\hom{R}( \quad  ,R)$. 

The projective stabilization $\bar{\mod R}$ is defined as follows. 
\begin{itemize} 
\item Each object of $\underline{\mod R}$ is an object of ${mod~R}$. 
\item For $A,B\in {\mod R}$, a set of morphisms from $A$ to $B$ is 
$\hombar{R}(A,B) = \hom{R}(A,B)/{\cal P}(A,B)$ where 
${\cal P}(A,B):= \{ f\in \hom{R}(A,B) \mid \mbox{$f$ factors through some projective module} \}$.  
Each element is denoted as $\underline{f} = f~\mod {\cal P}(A,B)$. 
A morphism $f:A\to B$ in $\mod R$ is called {\em a stable isomorphism} if $\bar{f}$ is an isomorphism in 
$\bar{\mod R}$ and we write $A\stcong B$.    
\end{itemize} 
For an $R$-module $M$, define a transpose $\tr M$ of $M$ 
to be $\cok \delta ^* $ where $P\stackrel{\delta }{\A} Q\A M\A 0$ is a projective 
presentation of $M$. The transpose of $M$ is uniquely determined as an object of $\bar{\mod R}$. 
If $\bar{f} \in \hombar{R}(M,N)$, then $f$ induces a map $\tr N \A \tr M$, which represents a morphism 
$\tr \bar{f} \in \hombar{R} (\tr N , \tr M )$. 
 
A kernel of projective cover of $M$ is called the first syzygy module of 
$M$ and denoted as $\syz{1}{R}{M}$. 
The first syzygy module of $M$ is uniquely determined as an object of 
$\bar{\mod ~R}$.  Inductively, we define 
$\syz{1}{R}{M}= \syz{1}{R}{\syz{n-1}{R}{M}}$. 
If $\bar{f} \in \hombar{R}(M,N)$, then $f$ induces a map 
$\syz{n}{R}{M} \A \syz{n}{R}{N}$, which represents a morphism 
$\bar{\syz{n}{R}{f}} \in \hombar{R} (\syz{n}{R}{M} , \syz{n}{R}{N} )$.

\begin{Lem} 
\label{p411}
On the commutative diagram with exact rows in $\mod R$
\[ \matrix{0&\A &A &\stackrel{f}{\A} &B &\stackrel{g}{\A} &C&\A &0\cr 
& &\mapdown{\alpha } & &\mapdown{\beta} & &\mapdown{\gamma}&&\cr
0&\A &A' &\stackrel{f'}{\A} &B' &\stackrel{g'}{\A} &C'&\A &0, \cr} \]
if $\beta$ and $\gamma$ are stable isomorphisms, so is $\alpha$. 
\end{Lem}

\noindent
{\bf proof.}  
We show that $\alpha$ is a stable isomorphism in the following case:
\begin{listing}{}
\item $\gamma$ is an isomorphism and $\beta$ is a stable isomorphism. 
\item $\beta$ and $\gamma$ are stable isomorphisms and 
$\gamma$ is an epimorphism. 
\item $\beta$ and $\gamma$ are stable isomorphisms. 
\end{listing}

\noindent 
1) Adding a projective cover $\rho _A: P_A \to A$ to the given diagram, 
we get the following: 
\[ \matrix{&&0&&0&&&&\cr &&\downarrow&&\downarrow&&&&\cr 
&&Q&=&Q&&&&\cr &&\mapdown{v}&&\mapdown{w}&&&&\cr
0&\A &A \oplus P_A &\stackrel{({{f0}\atop{01}})}{\A} &
B\oplus P_A &\stackrel{(g~0)}{\A} &C&\A &0\cr 
& &\mapdown{({\alpha ~ \rho _A }) } & &\mapdown{(\beta ~ f' \circ\rho _A )} & &\mapdown{\cong}&&\cr
0&\A &A' &\stackrel{f'}{\A} &B' &\stackrel{g'}{\A} &C'&\A &0 \cr
& &\mapdown{} &&\mapdown{} &&&& \cr
& &0 &&0 &&&& \cr} \] 
Since $(\beta ~ f' \circ\rho_A )$ is an epimorphism and a stable isomorphism at the same time, 
it is a split epimorphism and $Q$ is projective.  
In other words, $w$ is a split monomorphism. Hence $v$ is also 
a split monomorphism and 
$(\alpha  ~  \rho _A )$ a split epimorphism with a projective kernel. 
In particular, $\alpha$ is a stable isomorphism. 

\noindent 
2) Since $\gamma$ is a split epimorphism with a projective kernel, 
there exists $\gamma ' : C' \to C$ such that 
$\gamma \circ \gamma ' = \id _{C'}$. 
On the diagram 
\[ \matrix{0&\A &A &{\A} &B'' &{\A} &C'&\A &0\cr 
& &{\| } & &\mapdown{\beta '} & &\mapdown{\gamma '}&&\cr
0&\A &A &\stackrel{f}{\A} &B &\stackrel{g}{\A} &C &\A &0, \cr} \]
$\beta '$ is a stable isomorphism 
because 
$\ext{1}{R}(\gamma ' , A): \ext{1}{R}(C , A ) \to \ext{1}{R}(C' , A)$ 
is an isomorphism. 
Connection of the above and the given diagrams yields 
another diagram: 
\[ \matrix{0&\A &A &{\A} &B'' &{\A} &C'&\A &0\cr 
& &\mapdown{\alpha } & &\mapdown{\beta  \circ \beta '} & &
{\|}&&\cr
0&\A &A' &\stackrel{f'}{\A} &B' &\stackrel{g'}{\A} &C'&\A &0 \cr} \] 
Since $\beta  \circ \beta '$ is a stable isomorphism, 
we can apply 1). 

\noindent 
3) Adding a projective cover $\rho _{B'} : P_{B'} \to B'$, 
we get 
\[ \matrix{0&\A &A &\stackrel{f\choose 0}{\A} &B \oplus P_{B'} &
\stackrel{({{g0}\atop{01}})}{\A} &C \oplus P_{B'}&\A &0 \cr
& &\mapdown{\alpha } & &\mapdown{( \beta ~ \rho _{B' })} 
& &\mapdown{( \gamma  ~ g'\circ\rho _B )}&&\cr
0&\A &A' &\stackrel{f'}{\A} &B' &\stackrel{g'}{\A} &C'&\A &0. \cr } \]  
Since $\beta$ is a stable isomorphism and 
$( \gamma  ~ g'\circ\rho _B )$ is an epimorphic stable isomorphism, 
we can apply (2) to get the conclusion. 
\quad (q.e.d.) 

\vspace{2mm} 
   
Let $\ll$ be a full subcategory of $\hpy ({\mod R} )$ defined as 
\[ \ll = \{ F\Dot \in \hpy ( {\proj R} ) \mid 
\H{i}{F}=0~(i < 0), \quad \mathop{{\rm H} _j} 
( {F^*}_{\bullet} )=0~(j\ge 0 ) 
\}. \]    

\begin{Lem}
[ \cite{Okayama} Proposition 2.3, Proposition 2.4 ]
\label{p229, p230-1} 
~~
\begin{listing}{}
\item For $A \in \bar{\mod R}$, there exists ${F_A}\Dot \in \ll$ that satisfies 
\[ \H{0}{\tau _{\le 0} F_A} \stcong A. \]  Such an ${F_A}\Dot$ is uniquely 
determined by $A$ up to isomorphisms. 
We fix the notation ${F_A}\Dot$ and call this a standard resolution of $A$.
\item  For $\bar{f} \in \hombar{R}(A,B)$, there exists 
$f\Dot \in \hom{\hpy ( \mod R )}({F_A}\Dot ,{F_B}\Dot )$ that satisfies 
\[ \bar{\H{0}{\tau _{\le 0} f}}=\bar{f}. \]  Such an ${f}\Dot$ is uniquely 
determined by $\bar{f}$ up to isomorphisms, so we use the notation $f\Dot$ 
to describe a chain map with this property for given $\bar{f}$. 
\end{listing} 
\end{Lem}

\begin{Th}[ \cite{Okayama} Theorem 2.6]\label{p231-1}
The mapping $A \mapsto {F_A}\Dot$ gives a functor from $\bar{\mod R}$ to 
$\hpy (\mod R )$, and this gives a category equivalence 
between $\bar{\mod R}$ and $\ll$. 
\end{Th} 

For $f\in \hom{R}(A,B)$, there exists a triangle 
\begin{equation}\label{triangle}
  C(f) ^{\bullet -1} \stackrel{n\Dot}{\A} {F_A}\Dot  \stackrel{f\Dot}{\A} 
  {F_B}\Dot  \stackrel{c\Dot}{\A} C(f)\Dot . 
\end{equation} 
In general, $C\Dot$ does not belong to $\ll$ any more but it satisfies the following: 
\[ \H{i}{C}=0~(i < -1), \quad \mathop{{\rm H}_j} ( {C^*}_{\bullet} )=0~(j > -1 ) . \]

\begin{Def-Lem}[\cite{Okayama}, Definition and Lemma 3.1] As objects of $\bar{\mod R}$, 
$\bar{\ker}~\bar{f} := \H{-1}{\tau _{\le -1}C(f)}$ and 
$\bar{\cok}~\bar{f} := \H{0}{\tau _{\le 0}C}$ are uniquely determined by $\bar{f}$, 
up to isomorphisms. 
We call these the pseudo-kernel and the pseudo-cokernel of $\bar{f}$.  
And we have 
\[  \bar{\cok}~\bar{f} = {\tr}~\bar{\ker}~{\tr} \bar{f}.  \] 
\end{Def-Lem}

\noindent
{\bf Notations.} The pseudo-kernel or the pseudo-cokernel is determined as a stable equivalence class 
of modules, not as an isomorphic class of $\mod R$. Although, for the simplicity, a module $M$ with 
$M \stcong \bar{\ker}\bar{f}$ or $M \stcong \bar{\cok}\bar{f}$ is denoted by $\bar{\ker}\bar{f}$ or $\bar{\cok}\bar{f}$ 
respectively.  

\vspace{2mm} 

Chain maps $n\Dot$ and $c\Dot$ in (\ref{triangle}) induce 
$c_f : \bar{\cok}~\bar{f} \to A$ and $n_f : B \to \bar{\ker}~\bar{f} $. These maps have the properties 
$\bar{f}\circ \bar{n_f} = \bar{0}$ and $\bar{c_f}\circ \bar{f} = \bar{0}$. 
The following lemma shows why $\bar{\ker}~\bar{f}$ 
and $\bar{\cok} \bar{f}$ are called 
the pseudo-kernel and the pseudo-cokernel respectively. 

\begin{Lem}[\cite{Okayama} Lemma~ 3.3, Lemma~3.5] 
Let $f:A\to B$ be a homomorphism of $R$-modules. 
\begin{listing}{}
\item If $\bar{x} \in \hombar{R}(X,A)$ satisfies 
$\bar{f} \circ \bar{x} = \bar{0}$, there exists 
$\bar{h_x} \in \hombar{R}(X,\bar{\ker}~\bar{f})$ 
such that $\bar{x} = \bar{n_f} \circ \bar{h_x}$. 
\item If $\bar{y} \in \hombar{R}(B,Y)$ satisfies 
$\bar{y} \circ \bar{f} = \bar{0}$, there exists 
$\bar{e_y} \in \hombar{R}(\bar{\cok}\bar{f},Y)$ 
such that $\bar{y} = \bar{e_y}  \circ \bar{c_f}$. 
\end{listing}
\end{Lem}

 From (\ref{triangle}), we have an exact sequence 
\[ 0\to \bar{\ker}~\bar{f} \to A\oplus P \stackrel{(f~\rho )}{\to } B\to 0 \] 
with some projective module $P$. 
This characterizes the pseudo-kernel. 
  
\begin{Lem}
\label{p336} 
For a given $f \in \hom{R}(A,B)$, suppose both 
$A\oplus P \stackrel{\left( f~p \right) }{\A} B$ and 
$A\oplus P' \stackrel{\left( f~p' \right) }{\A} B$ are epimorphisms 
with projective modules $P$ and $P'$. 
Then there are stable isomorphisms
$\lambda : A \oplus P \to A \oplus P'$ and 
$\kappa : \ker {\left( f~p \right) } \A \ker {\left( f~p \right) }$ 
that make the following diagram commutative:
\[ \matrix{ 
0&\A & \ker {\left( f~p \right) } &\A &A\oplus P & 
\stackrel{\left( f~p \right) }{\A}& B & \A &0 \cr 
&& \mapdown{\kappa} &&\mapdown{\lambda}& 
&\| & &\cr 
0&\A & \ker {\left( f~p' \right) } &\A &A\oplus P' & 
\stackrel{\left( f~p' \right) }{\A}& B & \A &0 \cr } \]

\end{Lem}

\noindent {\bf proof.} Set $\pi : B\A \cok f$ as a canonical map. 
Both $\pi \circ p$ and $\pi \circ p'$ are 
projective covers of $\cok f$, hence there 
exists $l:P\A P'$ such that $\pi \circ p = \pi \circ p' \circ l$. 
Then $\Im {\left( p' \circ l -p \right)} \subset \Im f$, so we get 
$h:P \A \Im f$ as $h$ coincides with $p'\circ l -p$. 
Via $A\A \Im f$ which is surjective, $h$ can be lifted to a map 
$j:P\A A$. This shows the equation 
$f\circ j + p' \circ l =p$. 
The map $\lambda = \pmatrix{1&j \cr 0& l \cr}: A\oplus P \to A \oplus P'$ 
yields the desired diagram . 
Obviously $\lambda$ is a stable isomorphism, which 
implies that
$\kappa$ is a stable isomorphism from Lemma~\ref{p411}. 
(q.e.d.)

\vspace{2mm}
\begin{Lem}[\cite{Okayama} Lemma 3.6] \label{p138} 
~~\begin{listing}{}
\item  There is an exact sequence 
\[ \ex{\ker f }{\bar{\ker} ~f}{\syz{1}{R}{\cok ~f}}. \] 
 
\item There is an exact sequence 
\[ \ex{L}{\bar{\cok}~\bar{f}}{ \cok~{f}} \] 
such that $\syz{1}{R}{L}$ is the surjective image of $\ker{f}$.
\end{listing}
\end{Lem}

\begin{Lem}
\label{p362}
The following holds for $f \in \hom{R} (A,B)$.  
\begin{listing}{} 
\item $\bar{\ker}\bar{f}$ is projective if and only if 
$f\Dot$ can be taken as $ f^i$ are isomorphisms for $i \le -1$. 
\item If $\bar{\ker} \bar{f}$ is projective, 
then $\syz{1}{R}{f}$ is a stable isomorphism. 
\end{listing} \end{Lem} 

\noindent 
{\bf proof.} 
1) The "if" part is obvious. 
First notice that $\bar{\ker}\bar{f}$ is projective if and only if 
we can choose $C(f)\Dot$ such that 
$C(f) ^i =0 \quad ( i \le -2 )$ 
as an element of $\hpy (\mod ~R )$.  
There is an exact sequence 
\[ 0\to C(f) ^{\bullet -1} \to E\Dot \stackrel{{f'} \Dot}{\to} 
F\Dot \to 0 \] 
where $E\Dot \cong {F_A}\Dot$, $F\Dot \cong {F_B}\Dot$ 
in $\hpy ({\mod ~R})$ and via these isomorphisms, 
$f\Dot$ is isomorphic to ${f'}\Dot$. 
Easily we see that $C(f)^i =0$ if and only if 
${f'}^i$ is an isomorphism in $\mod~R$. 

2) The triangle 
\[ C(f)^{\bullet -1} \stackrel{n\Dot}{\A} {F_A}\Dot 
\stackrel{f\Dot}{\A} {F_B}\Dot \stackrel{c\Dot}{\A}  {C(f)}\Dot \]
induces an exact sequence of complexes  
\[ \ex{ {F_A}\Dot}{  {\tilde{F_B}}\Dot}
{  {C(f)}\Dot}  \] which again induces 
\[ \ex{\tau _{\le -1}  {F_A}\Dot}{\tau _{\le -1}  {\tilde{F_B}}\Dot}
{\tau _{\le -1}  {C(f)}\Dot}  \] 
where $ {\tilde{F_B}}\Dot$ is a direct sum of ${F_B}\Dot$ 
and a trivial complex. 
Hence the exact sequence of homology groups is 
\[ 0\to \syz{1}{R}{A} \stackrel{\syz{1}{R}{f} \choose  \star}{\to}
\syz{1}{R}{B}\oplus P \A \bar{\ker}\bar{f}  \to 0 \] 
with projective modules P. 
From the assumption, 
$\syz{1}{R}{f} \choose  \star$ is a split monomorphism, 
in particular, $\syz{1}{R}{f}$ is a stable isomorphism. 
\quad (q.e.d.) 

\begin{Cor}
\label{p369}
Suppose $R$ is local. If a morphism $f \in \mathrm{End}_{R} (A)$ 
satisfies that $\bar{\ker} \bar{f}$ is projective, 
then ${f}$ is a stable isomorphism. 
\end{Cor} 

\noindent
{\bf proof.} 
We may assume 
$C(f) ^i =0 \quad ( i \le -2 )$. 
Therefore the complex ${C(f)^*_\bullet}$ has no cohomology 
except for ${\mathrm H}_{-1} ( {C(f)}^*_\bullet )$. 
The triangle 
\[ {C(f)}_\bullet ^*  \to {F_A}^* _\bullet \to {F_A}^* _\bullet \to 
{C(f)}_{\bullet -1}^* \] 
induces an exact sequence 
\[ \ex{\mathrm{H}_{-1} (C(f)^*_\bullet )}{\mathrm{H}_{-1} ({F_A}^* _{\bullet })}
{\mathrm{H}_{-1} ({F_A}^* _{\bullet})}.  \] 
Since $R$ is local, a surjective endomorphism on a finite module 
is always an automorphism. Thus 
we get ${\mathrm{H}_{-1} (C(f)^* _\bullet )}=0$. It follows that 
$C(f)^*_{\bullet}$ is an exact sequence of projective modules, 
equivalently, $f$ is a stable isomorphism.  
\quad (q.e.d.)

\section{Representation by monomorphisms and perfect exact sequences}

\begin{Def}
A morphism $f:A \A B$ in $\mod R$ is said to be {\em 
represented by monomorphisms} (rbm for short) if 
some monomorphism $f':A' \A B'$ in $\mod R$ is projective-stably equivalent to $f$, 
that is, there exist stable isomorphisms $\alpha :A\A A'$  and 
$\beta :B \A B'$ such that $\bar{\beta \circ f} = \bar{ f' \circ \alpha}$. 
\end{Def}

Each morphism is not always rbm. 

\begin{Ex} 
Let $R$ be a ring of dimension $n \ge 3$, 
$N$ an $R$-module with $\pd N =n$, and 
$\varphi _N : N \to N^{**}$ the natural map. 
Then any map $ N\oplus P \to N^{**} \oplus Q$ of the form 
$\pmatrix{\varphi _N& \star\cr \star&\star\cr }$ 
with projective modules $P$ and $Q$, 
is never be monomorphic. If otherwise, $N\oplus P$ is 
a submodule of a projective module; this is a contradiction 
because $N$ has a maximal projective dimension. 
\end{Ex}

\vspace{2mm} 

It was Aulander and Bridger who first defined and studied 
"represented by monomorphisms" property.

\begin{Th}[Auslander-Bridger]\label{Auslander-Bridger}
The following are equivalent for a morphism 
$f: A\A B$ in $\mod R$.  
\begin{listing}{}
\item 
There exists a monomorphism 
$f': A\A B\oplus P$ with a projective module $P$ 
such that $f = s\circ f'$ via some split epimorphism $s : B\oplus P \A B$. 
\item
There exists a monomorphism 
$f': A\A B\oplus P$ with a projective module $P$ 
such that $f = s\circ f'$ via some split epimorphism $s : B\oplus P \A B$, 
and ${f'}^*$ is an epimorphism. 
\item 
$\hombar{R} (B,I) \A \hombar{R} (A,I)$ is surjective if 
$I$ is an injective module.
\end{listing}
\end{Th}

Auslander and Bridger's original definition of 
"represented by monomorphisms" condition is 1) of 
Theorem~\ref{Auslander-Bridger}. Seemingly this is different from 
our definition. But we show that two conditions are equivalent. 

\begin{Lem}
\label{p320-2}
For a morphism $f: A\A B$ in $\mod R$, $f$ is rbm 
if and only if there exists a monomorphism 
$f': A\A B\oplus P$ with a projective module $P$ 
such that $f = s\circ f'$ via some split epimorphism $s : B\oplus P \A B$. 
\end{Lem}

\noindent {\bf proof.} The " if" part is clear. We shall show "only if" part. 
Suppose there exists a monomorphism $f': A'\A B'$, 
stable isomorphisms $\alpha  : A \A A'$  and 
$\beta  : B \A B'$ such that 
$\bar{\beta \circ  f }=\bar{ f' \circ \alpha} $.
We first take projective covers $\rho _A : P_A \A A$ 
and $\rho _B : P_B \A B$ such that the induced map 
$f_P : P_A \to P_B$ by $f$ is a monomorphism. 
Since $\alpha$ is a stable isomorphism, there exists 
a morphism $\alpha ' : A' \to A$ such that 
$\bar{\alpha  \circ \alpha '} =\bar{ \mathrm{id}_{A'}}$ and   
$\bar{\alpha '  \circ \alpha} =\bar{ \mathrm{id}_{A}}$. 

From the last equation there exists a morphism 
$s_A : A \to P_A$ such that 
$\alpha ' \circ \alpha + \rho _A \circ s_A = \mathrm{id}_A$, 
equivalently $(\alpha '  \rho _{A} ) \circ 
{{\alpha}\choose{s_A}} = \mathrm{id}_A$. 
In particular, $(\alpha ' ~ \rho _{A} ): A' \oplus P_A \to A$ is a 
split epimorphism and ${{\alpha}\choose{s_A}} : A \to A' \oplus P_A$ 
is a split monomorphism. 
Similarly we get morphisms $\beta ' :B' \to B$, $s_B: B \to P_B$ and 
$s_{B'} : B' \to P_{B'}$ such that 
$(\beta ' ~  \rho _{B} ) \circ {{\beta}\choose{s_B}}
= \mathrm{id}_B$ and 
$(\beta ~  \rho _{B'} ) \circ {{\beta '}\choose{s_{B'}}}
= \mathrm{id}_{B'}$. 

Given equation $\bar {\beta \circ f} = \bar{f' \circ \alpha }$ 
induces 
$\bar{ \beta ' \circ  (  {\beta \circ f}  ) \circ  \alpha ' } 
= \bar{ \beta ' \circ  ( {f' \circ \alpha } ) \circ  \alpha ' }$, 
that is,  $\bar{f \circ \alpha '} = \bar{\beta ' \circ f'}$. 
Hence there exists a homomorphism 
$t: A' \to P_B$ such that 
$f\circ \alpha ' - \beta ' \circ f' = \rho _B \circ t$. 

Now we get a commutative diagram 
\[ \matrix{ A' \oplus P_A & 
\stackrel{\left( {{f'}\atop t} ~ {0\atop f_P} \right)}{\longrightarrow} & 
B' \oplus P_B \cr 
\mapdown{\left( \alpha ' ~ \rho _A \right) } && 
\mapdown{\left( \beta ' ~ \rho _B \right) }\cr 
A&\stackrel{f}{\to} & B.\cr }  \] 

Since  the composite of maps 
$\pmatrix{\beta ' &0\cr s_{B'} &0\cr 0& \mathrm{id}_{P_B} \cr} 
: B' \oplus P_B  \to B \oplus P_{B'} \oplus P_B$ and 
$\pmatrix{\beta &\rho _{B'}&0\cr 
0&0&\mathrm{id}_{P_B}\cr } : B \oplus P_{B'} \oplus P_B 
\to  B' \oplus P_B$ is equal to $\mathrm{id}_{B' \oplus P_B}$, 
$ \pmatrix{\beta '  & \rho _{B}\cr} \circ 
\pmatrix{f' & 0\cr t& f_P \cr} =  \pmatrix{\beta '  & \rho _{B}\cr} \circ
\mathrm{id}_{B' \oplus P_B} \circ  
\pmatrix{f' & 0\cr t& f_P \cr} = 
 \pmatrix{\beta '  & \rho _{B}\cr} \circ
 \pmatrix{\beta &\rho _{B'}&0\cr 
0&0&\mathrm{id}_{P_B}\cr } \circ 
\pmatrix{\beta ' &0\cr s_{B'} &0\cr 0& \mathrm{id}_{P_B} \cr} 
\circ \pmatrix{f' & 0\cr t& f_P \cr}$, 
and the following diagram commutes:
\[ \matrix{ A' \oplus P_A & 
\stackrel{f^{''}}{\longrightarrow} & 
B\oplus P_{B'} \oplus P_B \cr 
\mapdown{\left( \alpha ' ~ \rho _A \right) } && 
\mapdown{\rho  }\cr 
A&\stackrel{f}{\to} & B\cr }.  \] 
where
 \[  f'' = \pmatrix{\beta ' &0\cr s_{B'} &0\cr 0& \mathrm{id}_{P_B} \cr} 
\circ \pmatrix{f' & 0\cr t& f_P \cr}, \] and 
\[ \rho =  \pmatrix{\beta ' & \rho _{B}\cr} \circ
\pmatrix{\beta &\rho _{B'}&0\cr 
0&0&\mathrm{id}_{P_B}\cr }.   \] 
It is easy to see that $f''$ is a monomorphism and $\rho$ is a split 
epimorphism. 

Finally putting $f''' = f'' \circ \pmatrix{\alpha \cr s_A \cr}$ which is  
a monomorphism, we have 
\[ f = f \circ \mathrm{id}_A = f \circ \left( \alpha ' ~ \rho _A \right) 
\circ  \pmatrix{\alpha \cr s_A \cr} = \rho \circ f'' \circ  
\pmatrix{\alpha \cr s_A \cr} = \rho \circ f''' .  \] \quad 
(q.e.d.) 

\vspace{3mm}

The most remarkable point in Auslander-Bridger's Theorem is that 
being rbm is equivalent to 
being represented by "perfect monomorphisms" whose $R$-dual 
is an epimorphism. 

\begin{Def}
An exact sequence $\ex{A}{B}{C}$ of $R$-modules is called a perfect exact sequence or to be 
perfectly exact if its $R$-dual 
$\ex{\hom{R} (C,R)}{\hom{R} (B,R)}{\hom{R} (A,R)}$ is also exact. 
\end{Def}

\begin{Prop}[ \cite{Okayama} Lemma ~2.7]
\label{p366} 
The following are equivalent for an exact sequence 
\[ \theta : 0\A A \stackrel{f}{\A} B \stackrel{g}{\A} C \A 0.  \]
\begin{listing}{}
\item  $\theta $ is perfectly exact. 
\item  $0\to {F_A}\Dot  \stackrel{f\Dot}{\to} {{F_B}\Dot } \stackrel{g\Dot}{\to} {{F_C}\Dot }\to 0$@is exact.
\item ${F_C}^{\bullet -1} \A {F_A}\Dot \stackrel{f\Dot}{\to} 
{F_B}\Dot \stackrel{g\Dot}{\to} {F_C}\Dot$ is a distinguished triangle 
in $\hpy{ (\mod~R )}$. 
\item ${F_A}\Dot \cong C(g )^{\bullet -1} $ in $\hpy ({\proj R} )$. 
\end{listing}
If these conditions are satisfied, we have the following.
\begin{listing}{}
\setcounter{list}{4}
\item $ C \stcong \bar{\cok}{\bar{f}}$. 
\item ${F_C}\Dot \cong C(f)\Dot $ in $\hpy ({\proj R} )$. 
\end{listing}
\end{Prop}

\noindent
{\bf proof.} 

In  \cite{Okayama} Lemma ~2.7, we see the equivalence between 1) and 2). 

The implication 3) $\Rightarrow$ 2) is obvious. 
 
For the rest of the proof, consider the following diagram:
\[ \matrix{ 0&\A &{F_A}\Dot &{\A} & {\mathop Cyl}(f)\Dot &\A &{C(f)}\Dot & \A&0 \cr  
&&\|& &\mapdown{{\beta '} \Dot} & & \mapdown{{\gamma '}\Dot} &&\cr
0&\A &{F_A}\Dot &\stackrel{f\Dot}{\A} & {F_B}\Dot &\stackrel{g\Dot}{\A} &{F_C}\Dot & \A&0 \cr  
& &\mapdown{\alpha \Dot}& &\mapdown{{\beta ''}\Dot} &&\mapdown{{\gamma ''}\Dot} && \cr  
0&\A &{C(g)}^{\bullet -1} &\A & \widetilde{F_B} \Dot &\A &{\mathop Cyl}(g)\Dot & \A&0 \cr } \]
The top-row and the bottom-row are exact. Chain maps ${\beta '} \Dot$, 
${\beta ''} \Dot$ and ${\gamma ''}\Dot$ are 
isomorphisms up to homotopy. 

2) $\Rightarrow$ 3),4),5) and 6). 
If the middle row is also exact,  then 
$C(\gamma ')\Dot \cong C(\beta ')\Dot$, which are trivial complexes, 
hence ${\gamma '}\Dot$ is an isomorphism. 
Now ${\beta '' }\Dot \circ {\beta '}\Dot$ and 
 ${\gamma '' }\Dot \circ {\gamma '}\Dot$ are isomorphisms, 
$C(\alpha )\Dot =0$ follows from the exact sequence 
\[ \ex{ C(\alpha )\Dot}{{C({\beta '' } \circ {\beta '})}\Dot}
{{C({\gamma '' } \circ {\gamma '})}\Dot}.  \] 
 
4) $\Rightarrow$ 3). On the above diagram 
$\alpha ^i = \mathrm{id}$ for $i \le 0$, so 
${F_A}\Dot \cong C(g )^{\bullet -1} $ implies that $\alpha \Dot$ is 
an isomorphism. Therefore ${\gamma ''} \Dot\circ {\gamma '} \Dot$ is also an isomorphism. 
\quad (q.e.d.) 

\vspace{2mm}

If $R$ is local, all the conditions above are equivalent. 
We shall give the proof later at the end of this section.
 
\begin{Lem}\label{p366-2}
Let the sequence 
\[ \theta : 0\A A \stackrel{f}{\A} B \stackrel{g}{\A} C \A 0  \] 
be exact.
If $R$ is local, the conditions 1) - 4) in Proposition~\ref{p366} 
are equivalent to the conditions 5) and 6). 
\begin{listing}{}
\setcounter{list}{4}
\item $ C \stcong \bar{\cok}{\bar{f}}$. 
\item ${F_C}\Dot \cong C(f)\Dot $ in $\hpy ({\proj R} )$. 
\end{listing}
\end{Lem} 

 For a morphism $f:A \to B$, $A \oplus P_B \stackrel{(f~\rho _B )}{\longrightarrow} B$ is an epimorphism 
with a projective cover $\rho _B : P_B \to B$. Thus each morphism is represented by epimorphisms. 
And the choice of the representing epimorphism is unique up to direct sum of projective modules, 
as we have seen in Lemma~\ref{p336}. 

Unlikely, we already know an example of a morphism that is not rbm. 
And moreover, even if a given map is represented by a monomorphism, 
there would be another representing monomorphism.  
We see it in Example~\ref{p347}. 

However, uniqueness theorem is obtained in this way.
Due to Theorem ~\ref{Auslander-Bridger}, a morphism is rbm 
if and only if it is represented by a perfect monomorphism.  
And if this is the case, the representing perfect monomorphism is 
uniquely determined up to stable isomorphisms. 
These are the statements in Theorem~\ref{p357}, before which, 
we need some preparations. 

\vspace{2mm} 
 
For given exact sequence of modules 
$A\stackrel{f}{\to} B \stackrel{g}{\to} C$, 
we have a diagram of triangles 
 \begin{equation}\label{p355*}
  \matrix{ {F_A}\Dot& \stackrel{f\Dot}{\A}& {F_B}\Dot & 
\stackrel{}{\A}& {C(f)}\Dot &\stackrel{}{\A} & {F_A}\Dot \cr 
\mapdown{\alpha \Dot}&& \|&&\mapdown{\gamma \Dot}&&
\mapdown{\alpha ^{\bullet +1}} \cr 
C(g)^{\bullet -1}& \stackrel{}{\A}& {F_B}\Dot & 
\stackrel{g\Dot}{\A}& {F_C}\Dot &\stackrel{}{\A} & {C(g)}\Dot \cr } 
\end{equation} 
which induces a diagram with exact rows 
\begin{equation}\label{p355**} 
{\small \matrix{ 
0&\A & \H{-1}{C(f)}& \A&
 A &\stackrel{f\choose \epsilon }{\A}& B\oplus {F_A}^1 &
\stackrel{( c_f ~ \pi )}{\to}&\bar{\cok}\bar{f} & \A & 0  \cr 
& && &\mapdown{\alpha} & &\mapdown{\beta} &
&\mapdown{\gamma} & &  \cr
& &0& \A&\bar{\ker}\bar{g} &\to & B\oplus P_C &
\stackrel{(g~\rho _C )}{\to} &C & \A & 0  \cr } }
\end{equation}
We observe some facts below. 

\begin{Lem} 
\label{p356} 
With the notations above, the following holds. 
\begin{listing}{}
\item $\beta$ is a stable isomorphism. 
\item ${C(\alpha )}^{\bullet +1} \cong {C(\gamma )}\Dot$. 
\item $\alpha$ is the composite of natural maps 
$A \to \Im ~f = \ker ~g$ and $\ker ~g \to \bar{\ker}\bar{g}$. 
So if $f$ is injective and $g$ is surjective, 
then from Lemma~\ref{p138}, $\alpha$ is a stable isomorphism, 
$\tau _{\le -1} C( \alpha ) \Dot=0$, and 
$\tau _{\le -2} C( \gamma ) \Dot=0$. 
\item If $\H{-1}{C(f)} =0$, then the upper row of (\ref{p355**}) is 
the short exact sequence 
\[ \theta _f : 0\to A \stackrel{f\choose \epsilon }{\A} B\oplus {F_A}^1
\stackrel{( c_f ~ \pi )}{\to} \bar{\cok}\bar{f} \to 0 \] 
which is a perfect exact sequence. 
\end{listing}{}
\end{Lem}

\begin{Th}
\label{p357}
Let $f: A\A B$ be a morphism in $\mod R$. 
Then $f$ is rbm if and only if 
$\H{-1}{C(f)}$ vanishes. If this is the case, we have the following: 
\begin{listing}{}
\item   We have a perfect exact sequence 
\[ \theta _f : 0\to A \stackrel{f\choose \epsilon }{\A} B\oplus {F_A}^1
\stackrel{( c_f ~ \pi )}{\to} \bar{\cok}\bar{f} \to 0. \] 
\item For any exact sequence of the form 
\[ \sigma : 0\A A \stackrel{\left(  f\atop q \right) }{\A} B\oplus P' 
\stackrel{\left( g~p \right) }{\A} C \A 0  \]
with some projective module $P'$, 
there is a commutative diagram 
\[ \matrix{\theta _f:&0&\A&A&\stackrel{\left(  f\atop \epsilon \right) }{\A} &
B\oplus {F_A}^1&\stackrel{\left( c_f~\pi \right) }{\A}& 
\bar{\cok}\bar{f} &\A &0 \cr
&&&\mapdown{\tilde{\alpha}}& &\mapdown{\tilde{\beta}}& &
\mapdown{\tilde{\gamma}} && \cr
\sigma :&0&\A&A &\stackrel{\left(  f\atop q \right) }{\A} &B\oplus P'& 
\stackrel{\left( g~p \right) }{\A}& C &\A &0 \cr} \] 
where $\tilde{\alpha}$ and $\tilde{\beta}$ are stable isomorphisms. 
\item There is an exact sequence with some projective module $Q$ and $Q'$  
\[ 0\to {Q'}\to{\bar{\cok}\bar{f}} \oplus Q 
\stackrel{\left( \tilde{\gamma} ~\star \right )}{\to}{C} \to 0. \] 
\item If $\sigma$ is also perfectly exact, then 
$\sigma$ is isomorphic to $\theta _f$ up to 
direct sum of trivial complexes.    
\end{listing}
\end{Th} 

\noindent {\bf proof.}   
Suppose that $f$ is rbm; 
there is an exact sequence 
\[ \sigma :  0\A A \stackrel{\left(  f\atop q \right) }{\A} B\oplus P' 
\stackrel{\left( g~p \right) }{\A} C \A 0.  \]
The maps $\tilde{f} = {f\choose q}$ and 
$\tilde{g} = \left( g~p \right)$ produce 
the same diagram as  (\ref{p355*}) 
because we may consider 
$\tilde{f}\Dot = f\Dot$ and  $\tilde{g}\Dot = g\Dot$. 
Apply Lemma~\ref{p356} 3) to this sequence, and we get
$\tau _{\le -2}C(\gamma )\Dot =0$ as for 
$\gamma \Dot : {C(f)}\Dot \to {F_C}\Dot$. 
From the long exact sequence of homology groups 
$ \H{-2}{C(\gamma )} \to  \H{-1}{C(f )} \to  \H{-1}{F_C} $, 
we get $ \H{-1}{C(f )} =0$. 
Conversely, suppose that $ \H{-1}{C(f )} =0$. Then 
Lemma ~\ref{p356} 4) shows that $\theta _f$ is perfectly exact. 
Now it remains to prove 2) - 4) in the case $ \H{-1}{C(f )} =0$.

2) Applying the argument of Lemma~\ref{p356} to the sequence $\sigma$, 
we get a similar diagram as (\ref{p355**} ) :
\[ \matrix{ 
0&\A &
 A &\stackrel{\left( {f\atop q}\atop \epsilon \right)}{\A}& B\oplus P' \oplus {F_A}^1 &
\stackrel{}{\to}&\bar{\cok}\bar{f\choose q} & \A & 0  \cr 
& &\mapdown{\alpha} & &\mapdown{\beta} &
&\mapdown{\gamma} & &  \cr
0& \A&\bar{\ker}\bar{(g~p)} &\to & B\oplus P' \oplus P_C &
\stackrel{(g~p~\rho _C )}{\to} &C & \A & 0.  \cr } \] 
The upper row is a direct sum of $\theta _f$ and a trivial complex, 
and the lower row is that of $\sigma$ and a trivial complex. 
Hence we get a desired diagram. 
\[ \matrix{0&\A&A&\stackrel{\left(  f\atop \epsilon \right) }{\A} &
B\oplus {F_A}^1&\stackrel{\left( c_f~\pi \right) }{\A}& 
\bar{\cok}\bar{f} &\A &0 \cr
&&\mapdown{\tilde{\alpha}}& &\mapdown{\tilde{\beta}}& &
\mapdown{\tilde{\gamma}} && \cr
0&\A&A &\stackrel{\left(  f\atop q \right) }{\A} &B\oplus P& 
\stackrel{\left( g~p \right) }{\A}& C &\A &0. \cr} \] 
Notice that $\tilde{\beta}$ is a stable isomorphism. 
From Lemma~\ref{p356} 3), $\tilde{\alpha}$ is also a 
stable isomorphism.  

3) Consider the exact sequence of complex 
\[ 0\to {C(\gamma )}^{\bullet -1} \A \widetilde {C(f) }\Dot 
\stackrel{\gamma \Dot}{\to} {F_C}\Dot \to 0,  \]  
where $\widetilde {C(f) }\Dot \cong {C(f)}\Dot$ in $\hpy(proj R)$. 
Applying the truncation 
$\tau _{\le 0}$, we get 
\[ 0\to {(\tau _{\le -1} C(\gamma ) )}^{\bullet -1} 
\to \tau _{\le 0} {\widetilde{C(f)}}\Dot \to \tau _{\le 0} {F_C}\Dot \to 0,  \] 
which induces an exact sequence of homology 
$0\to Q' \to \bar{\cok}\bar{f} \oplus Q \to C\to 0$ 
with projective modules $Q' = C(\gamma ) ^{-1}$ and $Q$ 
from Lemma~\ref{p356} 3). 

4) Suppose $\sigma$ is perfect. From Proposition ~\ref{p366}, 
${F_C}^{\bullet -1} \A {F_A}\Dot \stackrel{f\Dot}{\to} 
{F_B}\Dot \stackrel{g\Dot}{\to} {F_C}\Dot$ is a 
distinguished triangle, and ${F_C}\Dot \cong {C(f)}\Dot$, hence the induced sequence 
$\sigma$ is isomorphic to $\theta _f$. 
\quad (q.e.d.) 
 
\vspace{2mm} 

\begin{Ex}
\label{p347} 
Let $k$ be a field and $R= k[[X,Y,Z]]/ (X^2 -YZ)$. Let $M$ be an $R$-module 
defined as $M=R/( XY, Y^2 , YZ )$.  
The minimal Cohen-Macaulay approximation  
$\theta _f : 0\to Y_k \stackrel{f}{\to} X_k \to k\to 0$ 
of $k$ is perfectly exact 
since $\ext{1}{R}(k, R) =0$. On the other hand, the minimal 
Cohen-Macaulay approximation (see \cite{AB} for definition) of $M$, 
$\sigma : 0\to Y_M \stackrel{g}{\to} X_M \to M\to 0$ is not perfect 
since $\ext{1}{R}(M, R) \neq 0$. The map $g$ is decomposed as 
$Y_M \cong Y_k \oplus R$,  $X_M \cong X_k \oplus R$ and 
$g = \left( \matrix{f&q\cr 0&Y\cr} \right)$. 
We easily check the statement of Theorem ~\ref{p357}; 
$k \stcong \bar{\cok}\bar{f}$ and 
there is an exact sequence 
\[ \ex{R}{k\oplus R}{M}  \] 
which clearly does not split.  
\end{Ex}

\noindent 
{\bf proof of Lemma~\ref{p366-2}.} 
5) $\Leftrightarrow $ 6). Since $f$ is injective, 
we get $\H{-1}{C(f)} =0$ from Theorem~\ref{p357}. 
In other words, $C(f)\Dot \cong F_{\bar{\cok}\bar{f}} \Dot$. 
In this situation, 5) and 6) are clearly equivalent. 

6) $\Rightarrow$ 4). Since the assumption implies $\H{-1}{C(f)} =0$, 
we get an exact sequence 
$0\to Q' \to C\oplus Q \stackrel{(\gamma ~\star )}{\to} C\to 0$ with 
projective modules $Q$ and $Q'$, applying Theorem ~\ref{p357} 3).  
Since $R$ is local, we can apply Corollary ~\ref{p369}, 
which shows that $\gamma$ is a stable isomorphism, 
equivalently $\gamma \Dot$ is an isomorphism, hence $\alpha \Dot$ is also 
an isomorphism.
 (q.e.d.) 

\vspace{2mm}

An exact sequence $\theta : \ex{A}{B}{C}$ is perfectly exact 
if $\ext{1}{R}(C,R) =0$. 
But as we see in the next example, the vanishing of $\ext{1}{R}(C,R)$ is not the sufficient condition 
for $\theta$ to be perfectly exact. 
@

\begin{Ex}
\label{p344}
Let $R$ be $R= k[[X,Y]] / (XY)$ and let $f$ be a map defined via 
projective resolutions as follows:
\[ \matrix{
\cdots& \A & R^2 &\stackrel{(X~Y)}{\longrightarrow} &R &\A &k \cr 
&  & \mapdown{(X~Y)} & &\mapdown{( {X\atop Y })} & &\mapdown{f} \cr 
0& \A & R&\stackrel{( {X\atop Y } )}{\longrightarrow} &R^2 &\A &Y^k \cr } \]
Then the sequence 
\[ 0\A X \A k \oplus R^2 \stackrel{(f~\rho _{Y^k} )}{\A} Y^k \A 0 \] 
with a projective cover $\rho _{Y^k}$ of $Y^k$ is a perfect exact sequence 
since $ \ext{1}{R}(k,R) \cong \ext{1}{R}(Y^k ,R)$. 
But $\ext{1}{R}(Y^k ,R) \neq 0$.
\end{Ex}

Notice that the dual of a perfect exact sequence is not always perfect as we see in the next example. 

\begin{Ex}
\label{p313} 
Let $(R,{\bf m}, k)$ be three-dimensional Gorenstein local ring. 
Set $L = \tr {k}$, $T_L = \syz{3}{R} \tr \syz{3}{R} \tr L$  
and $ \varphi _L : L \A T_L $ to be a natural map. 
Since $\H{i}{F_L} = 0 \quad ( i \le 1)$ and  
${F_{T_L}}\Dot$ is exact, $\H{-1}{C(\varphi _L )}=0$. 
Putting $N = \bar{\cok}\bar{\varphi _L}$, we have 
a perfectly exact sequence 
\[ \theta _{\varphi _L} : 0\A L \stackrel{{\varphi _L}\choose *}{\A} 
T_L \oplus P \A N\to 0 \] 
with some projective module $P$. It is easy to see that ${(N )}^*$ is free. 
Dualizing $\theta ^*$, we get an exact sequence 
\[ 0\A L^{**} \A {T_L}^{**} \A {(N)}^{**} \A \ext{1}{R}(L^* , R) \A 0.  \] 
But $\ext{1}{R} (L^* ,R) \cong  \ext{3}{R} (\tr L ,R) = 
\ext{3}{R}(k ,R) \neq 0 $. 
\end{Ex}

\begin{Rem}\label{p314}
Let $\theta : \ex{A}{B}{C}$ be a perfect exact sequence. 
Then $\theta ^*$ is also perfectly exact if and only if 
the induced map $ \H{2}{F_A} \A \H{2}{F_B}$ is a monomorphism. 
\end{Rem} 

\section{Representation by monomorphisms and torsionless modules.}

In the previous section, we see that a given map $f$ is represented by 
monomorphisms if and only if $\H {-1}{C(f)} =0$. If this is the case, 
$\bar{\ker} \bar{f} = \cok \dif{C(f)}{-2} $ is the first syzygy of 
$\bar{\cok} \bar{f} = \cok \dif{C(f)}{-1} $. 
So it is natural to ask the converse: Is a given map $f$ represented by 
monomorphisms if $\bar{\ker}\bar{f}$ is a first sygyzy? 
This section deals with the problem. 
As a conclusion, the answer is yes if 
the total ring of fractions $Q(R)$ of $R$ is Gorenstein. 
Notice that if $Q(R)$ is Gorenstein, then $Q(R)$ is 
Artinian as we see in Lemma~\ref{p419}. 
What is more, if $Q(R)$ is Gorenstein, instead of a pseudo-kernel,
we can use a (usual) kernel 
to describe rbm condition. 
We begin with seeing equivalent conditions for a module to be 
a first syzygy.

\begin{Def} 
An $R$-module $M$ is said to be torsionless
\footnote{In \cite{ABr}, Auslander and Bridger use the term "1-torsion free" for "torsionless". 
Usually a module $M$ is called torsion-free if the natural map $ M \to M \otimes Q(R)$ is injective. } 
if the natural map 
$\phi : M \to M^{**}$ is a monomorphism. 
\end{Def}

The next theorem is well known.  We use the proof in \cite{ABr} and \cite{EG}. 
\begin{Lem}
\label{p350}
The following are equivalent for an $R$-module $M$. 
\begin{listing}{}
\item $M$ is torsionless. 
\item $\ext{1}{R}(\tr{M} , R )=0$
\item $M$ is a first syzygy; there exists a monomorphism from $M$ to a projective module. 
\end{listing}{}
\end{Lem}

\noindent{proof.}
Let $\phi : M \to M^{**}$ be the natural map. 
The well known formula $\ker \phi \cong \ext{1}{R}(\tr{M},R)$ shows the equivalence between 1) and 2). 

2) $\Rightarrow$ 3). We may assume that $M$ is 
a submodule of a free module $R^l$. 
Let $ \cdots P^{-2} \stackrel{d^{-2}}{\to}  P^{-1} \stackrel{d^{-1}}{\to}  P^{0} \to \tr{M}$ be 
a free resolution of $\tr{M}$. Then 2) says 
${(P^0 ) }^* \stackrel{{(d^{-1})}^*}{\to}  {(P^{-1}) }^* \stackrel{{(d^{-2})}^*}{\to}  {(P^{-2 })}^*$ is exact. 
By definition, $M \cong \cok {(d^{-1})}^*$ which is isomorphic to $\syz{1}{R}{ \cok {d^{-2}}^*}$. 
 
3) $\Rightarrow$ 1).  We may assume that $M$ is a submodule of a free 
module $R^l$. 
Let $\left( \matrix{f_1 \cr f_2 \cr \vdots \cr f_l \cr} \right): M \to R^l$ be a monomorphisms.  
If $m \in M$ is not zero, $f(m)$ is not zero, so there exists some $i$ such that $0\neq f_i (m) = \phi (m) (f_i)$, 
which implies $\phi (m) \neq 0$. Thus $\ker ~\phi = (0)$. 
(q.e.d.) 

\begin{Lem}\label{p419} 
If $Q(R)$ is Cohen-Macaulay, then $Q(R)$ is of dimension zero. 
\end{Lem}

\noindent 
{\bf proof.} 
Suppose there exists a non-minimal prime ideal ${\bf q}$. 
Then ${\bf q} \not\in \ass Q(R)$, since $Q(R)$ has no 
embedded prime. This implies 
$\displaystyle{ {\bf q} \not\subset  \bigcup _{{\bf p} \in \ass Q(R)}  {\bf p}}$, 
hence ${\bf q}$ contains a non-zero-divisor which is a unit. 
(q.e.d.)

\begin{Lem}
\label{p385} 
Let $R$ be a Noetherian ring and $f: A \to B$ be a morphism in $\mod R$. 
Suppopse $Q(R)$ is Gorenstein. 
If $\bar{\ker}\bar{f}$ is projective, then $f$ is rbm.
\end{Lem}

\noindent 
{\bf proof.} 
The assumption says $\tau _{\le -2} C(f)\Dot =0$. From Theorem ~\ref{p357}, $f$ is rbm if and 
only if $\H{-1}{C(f)}=0$, which means that $\dif{C(f)}{-1}$ is injective. 
So we have only to show  $\ker \dif{C(f)}{-1} = {( \cok {(\dif{C(f)}{-1} )}^* )}^* =0$. 
A triangle 
${F_A}\Dot \stackrel{f\Dot}{\to} {F_B}\Dot \to {C(f)}\Dot \to 
{F_A}^{\bullet +1}$ induces an exact sequence 
\[ 0\to \mathrm{ H}_{-1} ( {C(f)}^*_{\bullet} ) \to \mathrm{ H}_{-1} ( {F_B}^*_{\bullet} ) \to  \mathrm{ H}_{-1} ( {F_A}^*_{\bullet} ) \to 0. \]
Note that $ \cok {(\dif{C(f)}{-1} )}^* \cong \mathrm{ H}_{-1} ( {C(f)}^* _\bullet )$ and  
$\mathrm{ H}_{-1} ( {F_B}^*_{\bullet} )  = \ext{1}{R}(B,R)$. 
Since $Q(R)$ is Gorenstein of dimension zero, 
${\ext{1}{R}(B,R) } \otimes Q(R) =0$. 
So if $\bf p$ is any associated prime ideal of $R$, 
${\ext{1}{R}(B,R) }_{\bf p} =0$. 
Hence  ${\cok {(\dif{C(f)}{-1} )}^*}_{\bf p} =0$, 
which implies  ${(\cok {(\dif{C(f)}{-1} )}^* )}^* =0$. 
\quad (q.e.d.) 

\vspace{2mm}

To solve our problem, the special kind of maps is a key. 
For $M \in \mod ~R$, consider a module $J^2 M = \tr \Omega ^1 _R \tr \Omega ^1 _R M$. Since $\tr J^2 M$ is 
a first sygyzy, we have $\ext {1}{R}(J^2 M , R ) =0$, which means $\mathrm{ H}_{-1} ({{F_{J^2 M}} ^*}_{\bullet} ) =0$ 
and $\tau _{\ge -2} {{F_{J^2 M}} ^*}_{\bullet} $is a projective resolution of 
$\tr \Omega ^1 _R M = \cok {( \dif{F_{J^2 M}}{-2} )}^* =
 \cok {( \dif{F_{ M}}{-2} )}^* $. 
The identity map on $\tr \Omega ^1 _R M $ induces a chain map 
$ {( F_{ M} ) }^*_{\bullet}  \to {( F_{J^2 M} ) }^*_{\bullet}$ 
and  
$ {\psi _M}\Dot :  { F_{J^2 M}  }\Dot \to  { F_{ M}  }\Dot$ 
subsequently. 
The maps $\psi _M ^i$ are identity maps for $i \le -1$, 
in other words, $\tau _{\le -2} C(\psi _M )\Dot =0$ and 
 $\bar{\ker}\bar{\psi _M}$ is projective. 
Thus we can apply  
the argument in Lemma~\ref{p385} for $f = \psi _M$; 
$\H{-1}{C(\psi _M )} \cong 
{( \mathrm{H}_{-1} {{C(\psi _M )^*_{\bullet}} })}^*$. 
Since $\mathrm{ H}_{-1} ({{F_{J^2 M}} ^*}_{\bullet} ) =0$, we have 
$\mathrm{ H}_{-1} (C(\psi _M ) ^* _{\bullet} ) \cong 
\mathrm{ H}_{-1}( (F_M ) ^* _{\bullet} ) 
\cong \ext{1}{R}(M,R) $. 
Therefore we have $\H{-1}{C(\psi _M)}  \cong 
{( \ext{1}{R}(M,R) ) }^* $. 
Now we get a result as follows: 

\begin{Cor}\label{p387} 
The map $\psi _M : J^2 M \to M$ is rbm 
if and only if an $R$-module $M$ has  ${\ext{1}{R}(M,R) }^* =0$.  
\end{Cor}

For given morphism of $R$-modules $f: A\to B$, 
adding a projective cover of $B$ to $f$, 
we get an exact sequence 
\[ 0\to \bar{\ker}\bar{f} \stackrel{{n_f}\choose \star}{\to} A\oplus P_B 
\stackrel{\left( f~\rho _B \right)}{\to} B \to 0.  \]
Due to Theorem~\ref{p357}, we have a perfect exact sequence 
$\theta _{n_f}$, because $n_f$ is rbm: 
\[ \matrix{
\theta _{n_f} :& 0&\to& \bar{\ker}\bar{f} &\stackrel{{n_f}\choose \diamond}{\to}& A\oplus F^1 _{\bar{\ker}\bar{f}} & \to& \bar{\cok}~\bar{n_f}&
 \to & 0\cr 
&&&\mapdown{\stcong} && \mapdown{\stcong} & 
&\mapdown{\omega _f}&  &\cr  
& 0&\to& \bar{\ker}\bar{f} &\stackrel{{n_f}\choose \star}{\to}& A\oplus P_B & 
\stackrel{\left( f~\rho _B \right)}{\to}& B& \to & 0\cr } \]

\begin{Lem}
\label{p391}
With notation as above, suppose 
${(\ext{1}{R}(\bar{\cok}\bar{f},R) )}^* =0$. 
Then the following conditions are equivalent. 
\begin{listing}{}
\item $f$ is rbm. 
\item $\bar{\ker}\bar{f}$ is torsionless and 
$\omega _f : \bar{\cok}~\bar{n_f} \to B$ is rbm. 
\end{listing}
\end{Lem} 

\noindent{\bf proof.} 
On the diagram of triangles 
\[ \matrix{
{F_{\bar{\ker}\bar{f}}}\Dot & \stackrel{{n_f}\Dot}{\to} & 
{F_A}\Dot &\to &{C(n_f)}\Dot & \to & {F_{\bar{\ker}\bar{f}}}^{\bullet +1}\cr 
\mapdown{\chi_f \Dot} & & 
\| & &\mapdown{\omega _f \Dot} && \mapdown{\chi_f ^{\bullet +1}}\cr 
{C(f)}^{\bullet -1} &\to & 
{F_A}\Dot &\stackrel{f\Dot}{\to} &{F_B}\Dot & \to & {C(f)}\Dot ,  \cr} \] 
we observe ${C(\chi_f )}^{\bullet +1} \cong {C(\omega _f )}\Dot$. 
We have 
$\H{-1}{C(n_f )} =0$ because $n_f$ is rbm. 
There is an exact sequence 
\[ \H{0}{F_{\bar{\ker}\bar{f}}} \to \H{-1}{C(f)} \to \H{-1}{C(\omega _f )}.  \] 

2) $\Rightarrow$ 1).  Since $\bar{\ker}\bar{f}$ is torsionless, 
$\H{0}{F_{\bar{\ker}\bar{f}}} =\ext{1}{R}(\tr \bar{\ker}\bar{f} ,R)=0$. 
And  $\H{-1}{C(\omega _f )}=0$ 
because $\omega _f$ is rbm. 
From the above exact sequence, we have $ \H{-1}{C(f)}=0$. 

1) $\Rightarrow$ 2). From the assumption, $ \H{-1}{C(f)}=0$ 
which implies $\bar{\ker}\bar{f} \stcong \syz{1}{R}{\bar{\cok}\bar{f}}$. 
We show that $\H{-1}{C(\omega _f )}=\H{0}{C(\chi_f )}$ vanishes. 
The equation $\bar{\ker}\bar{f} \stcong \syz{1}{R}{\bar{\cok}\bar{f}}$ 
implies ${F_{\bar{\ker}\bar{f}}}^{\bullet +1} \cong
 F\Dot_{J^2 ( \bar{\cok}\bar{f} )} $. On the other hand, 
$C(f)\Dot \cong {F_{\bar{\cok}\bar{f}}}\Dot$. 
Via these isomorphisms, $\chi_f^{\bullet +1}$ is regarded as 
$\psi \Dot _{\bar{\cok}\bar{f}}$. 
Hence from the proof of Corollary~\ref{p387}, 
$\H{0}{C(\chi_f )} \cong \H{-1}{C(\psi  _{\bar{\cok}\bar{f}} )}
\cong {( \ext {1}{R} ( \bar{\cok}\bar{f} , R) )}^* =0$. 
\quad(q.e.d.) 
  
\begin{Prop}
\label{p405}
Suppose $Q(R)$ is Gorenstein. 
A morphism $f$ of $R$-modules is rbm 
if and only if $\bar{\ker}\bar{f}$ is torsionless. 
\end{Prop}

\noindent {\bf proof.} 
In the previous section, we already have "only if" part. 
Apply Theorem~\ref{p357} to $n_f$ 
which is rbm, 
Theorem~\ref{p357} 3) says that  
$\bar{\ker}\bar{\omega _f}$ is projective. 
Therefore $\omega _f$ is rbm from 
Lemma~\ref{p385}. Since $Q(R)$ is Gorenstein, we can use 
Lemma~\ref{p391}, which completes the proof. \quad (q.e.d.)

\vspace{2mm} 
 
Now we go to the next stage; we are to state rbm 
condition in terms of normal kernel.

\begin{Lem}
\label{p326} 
Suppose $Q(R)$ is Gorenstein. 
Let the sequence of $R$-modules 
$0\A A\stackrel{f}{\A} B \stackrel{g}{\A} C \A 0 $ be exact. 
If $A$ and $C$ are torsionless, then so is $B$. 
\end{Lem} 

\noindent {\bf proof.} 
From the assumption, $A\cong \bar{\ker}\bar{g}$ is torsionless. 
Due to Proposition ~\ref{p405}, $g$ is rbm; there exists 
an exact sequence 
\[ \theta _g : 
0\A B \stackrel{g\choose q}{\A} C\oplus Q \A \bar{\cok}\bar{g} \A 0 \] 
with a projective module $Q$ and a map $q: B\A Q$. Since $C$ is 
a submodule of some projective module, so is $B$. 
(q.e.d.) 

\begin{Cor}
\label{p380} 
Suppose $Q(R)$ is Gorenstein. 
For a given morphism $f$, $\ker f$ is torsionless if and only if 
$\bar{\ker}\bar{f}$ is torsionless. 
\end{Cor}

\noindent {\bf proof.} 
From Lemma~\ref{p138}, there is an exact sequence 
$\ex{\ker f }{\bar{\ker} ~f}{\syz{1}{R}{\cok f}}$. 
So the "if" part is obvious, and the "only if" part comes from 
Lemma~\ref{p326}. 
\quad (q.e.d.)

\begin{Th}
\label{p284}
Suppose $Q(R)$ is Gorenstein. 
The following are equivalent for a morphism 
$f: A\A B$ in $\mod R$.  
\begin{listing}{}
\item \label{284-0}
$f$ is rbm. 
\item \label{284-9}
$\ker f$ is torsionless.
\item \label{284-6}
$\bar{\ker}\bar{f}$ is torsionless.
\item \label{284-4}
$\H{-1}{{C(f)}}=0$.
\item \label{284-5}
$\syz{1}{R}{\bar{\cok} \bar{f} } \stcong \bar{\ker}\bar{f} $ . 
\item \label{284-7}
There exists $f'$ such that $f' \stcong f$ and $\ker f'$ is torsionless. 
\item \label{284-8} 
For any $f'$ with $f' \stcong f$, $\ker f'$ is torsionless. 
\end{listing}\end{Th}

\noindent {\bf proof.} 
Implications $\ref{284-5})\Rightarrow \ref{284-6})$, 
$\ref{284-8})\Rightarrow \ref{284-9})$ and 
$\ref{284-8})\Rightarrow \ref{284-7})$ are obvious. 
We already showed $\ref{284-0})\Leftrightarrow \ref{284-4})$ 
in Theorem ~\ref{p357}, 
$\ref{284-0})\Leftrightarrow \ref{284-6})$ 
in Proposition ~\ref{p405}, and 
\ref{284-6})$\Leftrightarrow$ \ref{284-9}) 
in Corollary ~\ref{p380}. 
Implications $\ref{284-6})\Rightarrow \ref{284-8})$ and 
$\ref{284-7})\Rightarrow \ref{284-6})$ are obtained from 
"if" and "only if " part of Corollary~\ref{p380} respectively. 

$\ref{284-4})\Rightarrow \ref{284-5})$. It is clear since 
$\cok \dif{C(f)}{0}= \bar{\ker } \bar{f}$ and
$\cok \dif{C(f)}{-1}= \bar{\cok } \bar{f}$. 
\quad (q.e.d.)

\vspace{4mm}

The statement of Theorem~\ref{p284} does not hold 
for ring $R$ with $Q(R)$ non-Gorenstein. 

\begin{Cor}\label{p421}
The following are equivalent for a Noetherian ring $R$. 
\begin{listing}{}
\item $Q(R)$ is Gorenstein. 
\item Every morphism with torsionless kernel is 
rbm. 
\item ${\ext{1}{R}(M,R) }^* =0$ for each $M \in \mod R$. 
\end{listing}
\end{Cor} 

\noindent {\bf proof.} 

1) $\Rightarrow$ 2). It comes directly from Therem~\ref{p284}. 

For $M \in \mod R$, ${\ext{1}{R}(M,R) }^* =0$ means 
${\ext{1}{R_{\bf p}}(M_{\bf p},R_{\bf p}) } =0$ for 
every ${\bf p} \in \ass R$. 
Hence 3) says that ${R_{\bf p}}$ is 
a zero-dimensional Gorenstein for each ${\bf p} \in \ass R$, 
which is equivalent to 1) from Lemma ~\ref{p419}. 

2) $\Rightarrow$ 3). With no assumption, $\ker \psi _M$ is torsionless. 
Because $\bar{\ker}\bar{\psi _M}$ is projective and $\ker \psi _M$ is 
a submodule of $\bar{\ker}\bar{\psi _M}$ from Lemma ~\ref{p138} 1). 
So if 2) holds, $\psi _M$ is rbm and ${\ext{1}{R}(M,R) }^* =0$ 
from Corollary ~\ref{p387}.  (q.e.d.) 

\begin{Ex}
\label{p389} 
In the case $R = k[[X,Y,Z]]/(XY, X^2)$ with any field $k$, 
consider the map $\psi _k : k \to J^2 k$. 
We know $\bar{\ker}\bar{\psi _k}$ is projective. 
We shall show that $\psi _k$ is {\em not} rbm; 
that is, $\H{-1}{C(\psi _k )}$ does not vanish. 

As we have seen, 
${C(\psi _k )}^i = 0  \quad (i \le -2)$, 
$\H{-1}{C(\psi _k )} = {(\ext{1}{R} (k,R))}^*$. 
From a free resolution of $k$ 
\[ R^3 
\stackrel{\left( {Y\atop 0}{X\atop 0}{0\atop X} \right )}{\longrightarrow} 
R^2 \stackrel{\left( X~Y \right) }{\to} R \to k \to 0,  \]
we get a free resolution of $\ext{1}{R}(k,R)$  
\[ R^4 \stackrel{ \left( { {X~Y}\atop{0~0} }{ {0~0}\atop{X~Y} } \right) }
{\longrightarrow} R^2 \to \ext{1}{R}(k,R) \to 0.  \] 
We easily see that 
${( \ext{1}{R}(k,R) )} ^* \cong \H{-1}{C(\psi _k )}$ does not vanish. 
\end{Ex}

\vspace{2mm}

\noindent {\bf Acknowledgement.}  
I thank Kazuhiko Kurano who suggested that conditions for rbm 
should be given in terms of the kernel not only by the pseudokernel. 
I also thank Shiro Goto who told me that the assumption of 
Theorem~\ref{p284} is weakened and that Corollary~\ref{p421} 
holds.

\end{document}